# Embeddings of almost Hermitian manifolds in almost hyperHermitian those. Complex and hypercomplex numbers in differential geometry.

**Alexander A. Ermolitski**

*Cathedra of mathematics, MSHRC, st. Nezavisimosti 62,*
*Minsk, 220050, Belarus*
*E-mail: ermolitski@mail.by*

---

**Abstract:** Tubular neighborhoods play an important role in differential topology. We have applied these constructions to geometry of almost Hermitian manifolds. At first, we consider deformations of tensor structures on a normal tubular neighborhood of a submanifold in a Riemannian manifold. Further, an almost hyperHermitian structure has been constructed on the tangent bundle *TM* with help of the Riemannian connection of an almost Hermitian structure on a manifold *M* then, we consider an embedding of the almost Hermitian manifold *M* in the corresponding normal tubular neighborhood of the null section in the tangent bundle *TM* equipped with the deformed almost hyperHermitian structure of the special form.

As a result, we have obtained that any smooth manifold *M* of dimension *n* can be embedded as a totally geodesic submanifold in a Kaehlerian manifold of dimension 2*n* and in a hyperKaehlerian manifold of dimension 4*n*.

---

## 1. Deformations of tensor structures on a normal tubular neighborhood of a submanifold

**1°.** Let $(M',g')$ be a *k*–dimensional Riemannian manifold isometrically embedded in a *n*–dimensional Riemannian manifold $(M,g)$. The restriction of *g* to $M'$ coincides with *g*' and for any $p \in M'$.

$$T_p(M) = T_p(M') \oplus T_p(M')^{\perp}.$$

So, we obtain a vector bundle $M' \to T(M')^{\perp} : p \to T_p(M')^{\perp}$ over the submanifold $M'$. There exists a neighborhood $\widetilde{U}_0$ of the null section $O_{M'}$ in $T(M')^{\perp}$ such that the mapping

$$\pi \times \exp : v \to (\pi(v), \exp_{\pi(v)} v),\ v \in \widetilde{U}_0,$$

is a diffeomorphism of $\widetilde{U}_0$ onto an open subset $\widetilde{U} \subset M$. The subset $\widetilde{U}$ is called *a tubular neighborhood of the submanifold* $M'$ *in* $M$.

For any point $p \in M$ we can consider a set $\{\delta(p)\}$ of positive numbers such that the mapping $\exp_{U(\delta(p))}$ is defined and injective on $U(\delta(p)) \subset T_p(M)$. Let $\bar{\varepsilon}(p) = \sup\{\delta(p)\}$.

**Lemma, [6].** *The mapping* $M \to R_+ : p \to \bar{\varepsilon}(p)$ *is continuous on M.*

If we take the restriction of the function $\bar{\varepsilon}(p)$ on $\widetilde{U}$ then it is clear that there exists a continuous positive function $\varepsilon(p)$ on $M'$ such that for any $p \in M'$ open geodesic balls $B\left(p;\frac{\varepsilon(p)}{2}\right) \subset B(p;\varepsilon(p)) \subset \widetilde{U}$. For compact manifolds we can choose a constant function $\varepsilon(p) = \varepsilon > 0$. We denote $\widetilde{U}_p = \exp\left(\widetilde{U}_0 \cap T_p(M')^{\perp}\right)$, $D\left(p;\frac{\varepsilon(p)}{2}\right) = B\left(p;\frac{\varepsilon(p)}{2}\right) \cap \widetilde{U}_p$, $D(p;\varepsilon(p)) = B(p;\varepsilon(p)) \cap \widetilde{U}_p$. It is obvious that $\dim \widetilde{U}_p = \dim D(p;\varepsilon(p)) = n-k$. For any point $o \in M'$ we can consider such an orthonormal frame $\left(X_{1_0},\ldots,X_{n_0}\right)$ that $T_0(M') = L[X_{1_0},\ldots,X_{k_0}]$ and $T_0(M')^{\perp} = L[X_{k+1_0},\ldots,X_{n_0}]$. There exist coordinates $x_1,\ldots, x_k$ in some neighborhood $\widetilde{V}_0 \subset M'$ of the point *o* that $\frac{\partial}{\partial x_i}\Big|_0 = X_{i_0}, i = \overline{1,k}$. We consider orthonormal vector fields $X_{k+1}, \ldots, X_n$ which are cross–sections of the vector bundle $p \to T_p(M')^{\perp}$ over $\widetilde{V}_0$ and the neighborhood $\widetilde{W}_0 = \bigcup_{p \in \widetilde{V}_0} \cup \widetilde{U}_p$. The basis $\{X_{k+1_p},\ldots,X_{n_p}\}$ defines the normal coordinates $x_{k+1}, \ldots, x_n$ on $\widetilde{U}_p$ **[8]**. For any point $x \in \widetilde{W}_0$ there exists such unique point $p \in \widetilde{V}_0$ that $x = \exp_p(t\xi)$, $\|\xi\| = 1$, $\xi \in T_p(M')^{\perp}$. A point $x \in \widetilde{W}_0$ has the coordinates $x_1, \ldots, x_k, x_{k+1}, \ldots, x_n$ where $x_1, \ldots, x_k$ are coordinates of the point *p* in $\widetilde{V}_0$ and $x_{k+1}, \ldots, x_n$ are normal coordinates of *x* in $\widetilde{U}_p$. We denote $X_i = \frac{\partial}{\partial x_i}, i = \overline{1,n}$, on $\widetilde{W}_0$. Thus, we can consider *tubular neighborhoods* $Tb\left(M';\frac{\varepsilon(p)}{2}\right) = \bigcup_{p \in M'} D\left(p;\frac{\varepsilon(p)}{2}\right)$ and $Tb(M';\varepsilon(p)) = \bigcup_{p \in M'} D(p;\varepsilon(p))$ of the submanifold $M'$.

**2°.** Let *K* be a smooth tensor field of type (*r*, *s*) on the manifold *M* and for $x \in \widetilde{W}_0$, let

$$K_x = \sum_{i_1,\ldots,i_r,j_1,\ldots,j_s} k^{i_1,\ldots,i_r}_{j_1,\ldots,j_s}(x) X_{i_{1x}} \otimes \ldots \otimes X_{i_{rx}} \otimes X_x^{j_1} \otimes \ldots \otimes X_x^{j_s},$$

where $\{X_x^1,\ldots,X_x^n\}$ is the dual basis of $T_x^*(M)$, $x = \exp_p(t\xi)$, $\|\xi\| = 1$, $\xi \in T_p(M')^{\perp}$. We define a tensor field $\overline{K}$ on *M* in the following way.

a) $x \in D\left(p;\dfrac{\varepsilon(p)}{2}\right)$, then

$$\overline{K}_x = \sum_{i_1,\dots,i_r,j_1,\dots,j_s} k^{i_1,\dots,i_r}_{j_1,\dots,j_s}(p) X_{i_{1_x}} \otimes \dots \otimes X_{i_{r_x}} \otimes X^{j_1}_x \otimes \dots \otimes X^{j_s}_x;$$

b) $x \in D(p;\varepsilon(p)) \setminus D\left(p;\dfrac{\varepsilon(p)}{2}\right)$, then

$$\overline{K}_x = \sum_{i_1,\dots,i_r,j_1,\dots,j_s} k^{i_1,\dots,i_r}_{j_1,\dots,j_s}\left(\exp_p((2t-\varepsilon(p))\xi)\right) X_{i_{1_x}} \otimes \dots \otimes X_{i_{r_x}} \otimes X^{j_1}_x \otimes \dots \otimes X^{j_s}_x;$$

c) $x \in M \setminus \bigcup_{M'} D(p;\varepsilon(p))$, then

$$\overline{K}_x = K_x.$$

It is easy to see the independence of the tensor field $\overline{K}$ on a choice of coordinates in $\widetilde{W}_0$ for every point $o \in M'$.

**Definition 1.** *The tensor field $\overline{K}$ is called a deformation of the tensor field K on the normal tubular neighborhood of a submanifold $M'$.*

**Remark.** *The obtained tensor field $\overline{K}$ is continuous but is not smooth on the boundaries of the normal tubular neighborhoods $Tb\left(M';\dfrac{\varepsilon(p)}{2}\right)$ and $Tb(M';\varepsilon(p))$, $\overline{K}$ is smooth in other points of the manifold M.*

**3°.** We consider a deformation $\overline{g}$ of the Riemannian metric $g$ on the normal tubular neighborhood $Tb(M';\varepsilon(p))$ of a submanifold $M'$. For $x \in \widetilde{W}_0$, $x = \exp_p(t\xi)$, $\|\xi\| = 1$, $\xi \in T_p(M')$, we define the Riemannian metric $\overline{g}$ by the following way.

a) $\overline{g}_p = g_p$ for any $p \in M'$;

b) $\overline{g}_x(X_i, X_j) = \overline{g}_{ij}(x) = \overline{g}_{ij}(p)$, where $X_i = \dfrac{\partial}{\partial x_i}$, $i = \overline{1,n}$, $X_j = \dfrac{\partial}{\partial x_j}$, $j = \overline{1,n}$, on $\widetilde{W}_0$, $x \in D\left(p;\dfrac{\varepsilon(p)}{2}\right)$;

c) $\overline{g}_x(X_i, X_j) = \overline{g}_{ij}(x) = \overline{g}_{ij}\left(\exp_p((2t-\varepsilon(p))\xi)\right)$, for any $x \in D(p;\varepsilon(p)) / D\left(p;\dfrac{\varepsilon(p)}{2}\right)$;

d) $\overline{g}_x = g_x$ for each point $x \in M \setminus \bigcup_{p \in M'} D(p;\varepsilon(p))$.

The independence of $\bar{g}$ on a choice of local coordinates follows and the correctly defined Riemannian metric $\bar{g}$ on $M$ has been obtained.

It is known from **[9]** that every autoparallel submanifold of $M$ is a totally geodesic submanifold and a submanifold $M'$ is autoparallel if and only if $\nabla_X Y \in T(M')$ for any $X, Y \in \chi(M')$, where $\nabla$ is the Riemannian connection of $g$.

**Theorem 1.** *Let $M'$ be a submanifold of a Riemannian manifold (M, g) and $\bar{g}$ be the deformation of g on the normal tubular neighborhood $Tb(M';\varepsilon(p))$ of $M'$ constructed above. Then $M'$ is a totally geodesic submanifold of* $\left(Tb\left(M';\frac{\varepsilon(p)}{2}\right),\bar{g}\right)$.

**Proof.** For any point $x \in D\left(p;\frac{\varepsilon(p)}{2}\right) \subset \widetilde{W}_0$ the functions $\bar{g}_{ij}(x) = g_{ij}(p)$ and $\frac{\partial \bar{g}_{ij}}{\partial x_l} = 0,\quad l = \overline{k+1,n}$ on $D\left(p;\frac{\varepsilon(p)}{2}\right)$ because the vector fields $X_l = \frac{\partial}{\partial x_l}$ are tangent to $D\left(p;\frac{\varepsilon(p)}{2}\right)$. By the formula of the Riemannian connection $\bar{\nabla}$ of the Riemannian metric $\bar{g}$, **[8]**, we obtain for $i,j = \overline{1,k},\quad l = \overline{k+1,n}$

(1.1) $$2\bar{g}_p\left(\bar{\nabla}_{X_i} X_j, X_l\right) = X_{i_p}\bar{g}\left(X_j, X_l\right) + X_{j_p}\bar{g}\left(X_i, X_l\right) - X_{l_p}\bar{g}\left(X_i, X_j\right) + \\ + \bar{g}_p\left(\left[X_i, X_j\right], X_l\right) + \bar{g}_p\left(\left[X_l, X_i\right], X_j\right) + \bar{g}_p\left(X_i, \left[X_l, X_j\right]\right) = -\frac{\partial \bar{g}_{ij}}{\partial x_l} = 0.$$

Here we use the fact that $\left[X_i, X_j\right] = \left[X_l, X_i\right] = \left[X_l, X_j\right] = 0$ and that $\bar{g}\left(X_j, X_l\right) = \bar{g}\left(X_i, X_l\right) = 0$ because $X_l \in T(M')^{\perp}$.

Thus, $\bar{\nabla}_{X_i} X_j \in T(M')$ and from the remarks above the theorem follows.

**QED.**

**Corollary 1.1.** *Let $\bar{R}$ be the Riemannian curvature tensor field of $\bar{\nabla}$. Then $\bar{R}$ vanishes on every $D\left(p;\frac{\varepsilon(p)}{2}\right)$ for $p \in M'$.*

**Proof**. From the formula (1.1) it is clear that $\bar{\nabla}_{X_l} X_m = 0$ for $l,m = \overline{k+1,n}$. The rest is obvious.

**QED.**

## 2. Almost hyperHermitian structures (ahHs) on tangent bundles

**0°.** Let ($M$, $g$) be a $n$–dimensional Riemannian manifold and $TM$ be its tangent bundle. For a Riemannian connection $\nabla$ we consider the connection map $K$ of $\nabla$ **[2], [6],** defined by the formula

(2.1) $\nabla_X Z = KZ_* X$,

where $Z$ is considered as a map from $M$ into $TM$ and the right side means a vector field on $M$ assigning to $p \in M$ the vector $KZ_* X_p \in M_p$.

If $U \in TM$, we denote by $H_U$ the kernel of $K_{|TM_U}$ and this $n$–dimensional subspace of $TM_U$ is called the horizontal subspace of $TM_U$.

Let $\pi$ denote the natural projection of $TM$ onto $M$, then $\pi_*$ is a $C^\infty$–map of $TTM$ onto $TM$. If $U \in TM$, we denote by $V_U$ the kernel of $\pi_{*|TM_U}$ and this $n$–dimension subspace of $TM_U$ is called the vertical subspace of $TM_U$ $(\dim TM_U = 2\dim M = 2n)$. The following maps are isomorphisms of corresponding vector spaces $(p = \pi(U))$

$$\pi_{*|TM_U} : H_U \to M_p, K_{|TM_U} : V_U \to M_p$$

and we have

$$TM_U = H_U \oplus V_U$$

If $X \in \chi(M)$, then there exists exactly one vector field on $TM$ called the «horizontal lift» (resp. «vertical lift») of $X$ and denoted by $\overline{X}^h \left(\overline{X}^v\right)$, such that for all $U \in TM$ :

(2.2) $\pi_* \overline{X}^h_U = X_{\pi(U)}, \quad K\overline{X}^h_U = 0_{\pi(U)}$,

(2.3) $\pi_* \overline{X}^v_U = 0_{\pi(U)}, \quad K\overline{X}^v_U = X_{\pi(U)}$,

Let $R$ be the curvature tensor field of $\nabla$, then following **[2]** we write

(2.4) $[\overline{X}^v, \overline{Y}^v]=0$,

(2.5) $[\overline{X}^h, \overline{Y}^v] = \left(\overline{\nabla_X Y}\right)^v$

(2.6) $\pi_*\left([\overline{X}^h, \overline{Y}^h]_U\right) = [X, Y]$,

(2.7) $K\left([\overline{X}^h, \overline{Y}^h]_U\right) = R(X, Y)U$.

For vector fields $\overline{X} = \overline{X}^h \oplus \overline{X}^v$ and $\overline{Y} = \overline{Y}^h \oplus \overline{Y}^v$ on $TM$ the natural Riemannian metric $\hat{g} = <,>$ is defined on $TM$ by the formula

(2.8) $<\overline{X}, \overline{Y}> = g\left(\pi_* \overline{X}, \pi_* \overline{Y}\right) + g\left(K\overline{X}, K\overline{Y}\right)$.

It is clear that the subspaces $H_U$ and $V_U$ are orthogonal with respect to $<,>$.

It is easy to verify that $\overline{X}^h_1, \overline{X}^h_2, ..., \overline{X}^h_n, \overline{X}^v_1, \overline{X}^v_2, ..., \overline{X}^v_n$ are orthonormal vector fields on $TM$ if $X_1, X_2, ..., X_n$ are those on $M$ i.e. $g\left(X_i, X_j\right) = \delta^i_j$.

**1°.** We define a tensor field $J_1$ on $TM$ by the equalities

(2.9) $J_1\overline{X}^h = \overline{X}^v, J_1\overline{X}^v = -\overline{X}^h, X \in \chi(M)$.

For $X \in \chi(M)$ we get

$$J_1^2\overline{X} = J_1\left(J_1\left(\overline{X}^h \oplus \overline{X}^v\right)\right) = J_1\left(-\overline{X}^h \oplus \overline{X}^v\right) = -\left(\overline{X}^h \oplus \overline{X}^v\right) = -I\overline{X}$$

and

$$J_1^2 = -I.$$

For $X, Y \in \chi(M)$ we obtain

$$< J_1\overline{X}, J_1\overline{Y} > = < -\overline{X}^h \oplus \overline{X}^v, -\overline{Y}^h \oplus \overline{Y}^v > = < -\overline{X}^h, -\overline{Y}^v > + < \overline{X}^v, \overline{Y}^v >,$$

$$< \overline{X}, \overline{Y} > = < \overline{X}^h \oplus \overline{X}^v, \overline{Y}^h \oplus \overline{Y}^v > = < \overline{X}^h, \overline{Y}^h > + < \overline{X}^v, \overline{Y}^v >$$

and it follows that $< J_1\overline{X}, J_1\overline{Y} > = < \overline{X}, \overline{Y} >, (TM, J_1, <,>)$ is an almost Hermitian manifold.

Further, we want to analyze the second fundamental tensor field $h^1$ of the pair ($J_1$, <,>) where $h^1$ is defined by (2.11), **[3].**

The Riemannian connection $\hat{\nabla}$ of the metric $\hat{g} = <,>$ on $TM$ is defined by the formula (see **[6]**)

(2.10) $$< \hat{\nabla}_{\overline{X}}\overline{Y}, \overline{Z} > = \frac{1}{2}\left(\overline{X} < \overline{Y}, \overline{Z} > + \overline{Y} < \overline{Z}, \overline{X} > - \overline{Z} < \overline{X}, \overline{Y} > + \right.$$
$$\left. + < \overline{Z}, [\overline{X}, \overline{Y}] > + < \overline{Y}, [\overline{Z}, \overline{X}] > + < \overline{X}, [\overline{Z}, \overline{Y}] > \right), X, Y, Z \in \chi(M).$$

For orthonormal vector fields $\overline{X}, \overline{Y}, \overline{Z}$ on $TM$ we obtain

(2.11) $$h^1_{\overline{XYZ}} = < h^1_{\overline{X}}\overline{Y}, \overline{Z} > = \frac{1}{2} < \hat{\nabla}_{\overline{X}}\overline{Y} + J_1\hat{\nabla}_{\overline{X}}J_1\overline{Y}, \overline{Z} > =$$
$$= \frac{1}{2}\left(< \hat{\nabla}_{\overline{X}}\overline{Y}, \overline{Z} > - < \hat{\nabla}_{\overline{X}}J_1\overline{Y}, J_1\overline{Z} >\right) =$$
$$= \frac{1}{4}(< [\overline{X}, \overline{Y}], \overline{Z} > + < [\overline{Z}, \overline{X}], \overline{Y} > + < [\overline{Z}, \overline{Y}], \overline{X} > -$$
$$- < [\overline{X}, J_1\overline{Y}], J_1\overline{Z} > - < [J_1\overline{Z}, \overline{X}], J_1\overline{Y} > - < [J_1\overline{Z}, J_1\overline{Y}], \overline{X} >).$$

Using (2.4) – (2.7) and (2.11) we consider the following cases for the tensor field $h^1$ assuming all the vector fields to be orthonormal.

1.1°) $$h^1_{\overline{X}^h\overline{Y}^h\overline{Z}^h} = \frac{1}{4}(< [\overline{X}^h, \overline{Y}^h], \overline{Z}^h > + < [\overline{Z}^h, \overline{X}^h], \overline{Y}^h > +$$
$$+ < [\overline{Z}^h, \overline{Y}^h], \overline{X}^h > - < [\overline{X}^h, J_1\overline{Y}^h], J_1\overline{Z}^h > - < [J_1\overline{Z}^h, \overline{X}^h], J_1\overline{Y}^h > -$$
$$- < [J_1\overline{Z}^h, J_1\overline{Y}^h], \overline{X}^h >) = \frac{1}{4}\left(g([X,Y],Z) + g([Z,X],Y) + g([Z,Y],X) - \right.$$
$$\left. - < [\overline{X}^h, \overline{Y}^v], \overline{Z}^v > - < [\overline{Z}^v, \overline{X}^h], \overline{Y}^v > - < [\overline{Z}^v, \overline{Y}^v], \overline{X}^h > \right) =$$
$$= \frac{1}{2}g(\nabla_X Y, Z) - \frac{1}{4}\left(g(\nabla_X Y, Z) - g(\nabla_X Z, Y)\right) =$$
$$= \frac{1}{2}\left(g(\nabla_X Y, Z) - g(\nabla_X Y, Z)\right) = 0.$$

2.1°) $$h^1_{\overline{X}^h\overline{Y}^h\overline{Z}^v}=\frac{1}{4}\big(<[\overline{X}^h,\overline{Y}^h],\overline{Z}^v>+<[\overline{Z}^v,\overline{X}^h],\overline{Y}^h>+$$

$$+<[\overline{Z}^v,\overline{Y}^h],\overline{X}^h>-<[\overline{X}^h,J_1\overline{Y}^h],J_1\overline{Z}^v>-<[J_1\overline{Z}^v,\overline{X}^h],J_1\overline{Y}^h>-$$

$$-<[J_1\overline{Z}^v,J_1\overline{Y}^h],\overline{X}^h>\big)=\frac{1}{4}\Big(g(R(X,Y)U,Z)+<[\overline{Z}^h,\overline{X}^h],\overline{Y}^v>\Big)=$$

$$=\frac{1}{4}\big(g(R(X,Y)U,Z)+g(R(Z,X)U,Y)\big)=$$

$$=-\frac{1}{4}\big(g(R(X,Y)Z,U)+g(R(Z,X)Y,U)\big).$$

By similar arguments we obtain

3.1°) $h^1_{\overline{X}^h\overline{Y}^v\overline{Z}^h}=-\frac{1}{4}\big(g(R(Z,X)Y,U)+g(R(X,Y)Z,U)\big).$

4.1°) $h^1_{\overline{X}^v\overline{Y}^h\overline{Z}^h}=-\frac{1}{4}\big(g(R(Z,Y)X,U)\big).$

5.1°) $h^1_{\overline{X}^v\overline{Y}^v\overline{Z}^v}=\frac{1}{4}\big(g(R(Z,Y)X,U)\big).$

6.1°) $h^1_{\overline{X}^v\overline{Y}^v\overline{Z}^h}=0.$

7.1°) $h^1_{\overline{X}^v\overline{Y}^h\overline{Z}^v}=0.$

8.1°) $h^1_{\overline{X}^h\overline{Y}^v\overline{Z}^v}=0.$

It is obvious that $(J_1,\hat{g})$ is a Kaehlerian structure if and only if $h^1=0$.

**2°.** Now assume additionally that we have an almost Hermitian structure $J$ on $(M, g)$. We define a tensor field $J_2$ on $TM$ by the equalities

(2.12) $J_2\overline{X}^h=\left(\overline{JX}\right)^h,\quad J_2\overline{X}^v=-\left(\overline{JX}\right)^v,\quad X\in\chi(M).$

For $X\in\chi(M)$ we get

$$J_2^2\overline{X}=J_2\left(J_2\left(\overline{X}^h\oplus\overline{X}^v\right)\right)=J_2\left(\left(\overline{JX}\right)^h\oplus-\left(\overline{JX}\right)^v\right)=-\left(\overline{X}^h\oplus\overline{X}^v\right)-I\overline{X}$$

and

$$J_2^2=-I.$$

For $X,Y\in\chi(M)$ we obtain

$$<J_2\overline{X},J_2\overline{Y}>=<\left(\overline{JX}\right)^h\oplus-\left(\overline{JX}\right)^v,\left(\overline{JY}\right)^h\oplus-\left(\overline{JY}\right)^v>=<\left(\overline{JX}\right)^h,\left(\overline{JY}\right)^h>+$$

$$+<\left(\overline{JX}\right)^v,\left(\overline{JY}\right)^v>=g(JX,JY)+g(JX,JY)=g(X,Y)+g(X,Y)=$$

$$=<\overline{X}^h,\overline{Y}^h>+<\overline{X}^v,\overline{Y}^v>=<\overline{X}^h\oplus\overline{X}^v,\overline{Y}^h\oplus\overline{Y}^v>=<\overline{X},\overline{Y}>.$$

Further, we obtain

$$J_1\left(J_2\overline{X}\right)=J_1\left(\left(\overline{JX}\right)^h\oplus-\left(\overline{JX}\right)^v\right)=\left(\overline{JX}\right)^h\oplus\left(\overline{JX}\right)^v,$$

$$J_2\left(J_1\overline{X}\right)=J_2\left(-\overline{X}^h\oplus\overline{X}^v\right)=-\left(\overline{JX}\right)^h\oplus-\left(\overline{JX}\right)^v.$$

Thus, we get $J_1J_2 = -J_2J_1 = J_3$ and ahHs $(J_1, J_2, J_3, <,>)$ on $TM$ has been constructed.

For orthonormal vector fields $\overline{X}, \overline{Y}, \overline{Z}$ on $TM$ we obtain

(2.13) $$h^2_{\overline{XYZ}} = < h^2_{\overline{X}}\overline{Y}, \overline{Z} > = \frac{1}{2} < \hat{\nabla}_{\overline{X}}\overline{Y} + J_2\hat{\nabla}_{\overline{X}}J_2\overline{Y}, \overline{Z} > =$$
$$= \frac{1}{2}\left(< \hat{\nabla}_{\overline{X}}\overline{Y}, \overline{Z} > - < \hat{\nabla}_{\overline{X}}J_2\overline{Y}, J_2\overline{Z} >\right) = \frac{1}{4}(<[\overline{X}, \overline{Y}], \overline{Z} > +$$
$$+ < [\overline{Z}, \overline{X}], \overline{Y} > + < [\overline{Z}, \overline{Y}], \overline{X} > - < [\overline{X}, J_2\overline{Y}], J_2\overline{Z} > -$$
$$- < [J_2\overline{Z}, \overline{X}], J_2\overline{Y} > - < [J_2\overline{Z}, J_2\overline{Y}], \overline{X} >).$$

Using (2.4) – (2.7) and (2.13) we consider the following cases for the tensor field $h^2$ assuming all the vector fields to be orthonormal.

1.2°) $$h^2_{\overline{X}^h\overline{Y}^h\overline{Z}^h} = \frac{1}{4}(<[\overline{X}^h, \overline{Y}^h], \overline{Z}^h > + < [\overline{Z}^h, \overline{X}^h], \overline{Y}^h > +$$
$$+ < [\overline{Z}^h, \overline{Y}^h], \overline{X}^h > - < [\overline{X}^h, J_2\overline{Y}^h], J_2\overline{Z}^h > - < [J_2\overline{Z}^h, \overline{X}^h], J_2\overline{Y}^h > -$$
$$- < [J_2\overline{Z}^h, J_2\overline{Y}^h], \overline{X}^h >) = \frac{1}{4}(g([X,Y],Z) + g([Z,X],Y) + g([Z,Y],X) -$$
$$- g([X,JY],JZ) - g([JZ,X],JY) - g([JZ,JY],X)) =$$
$$= \frac{1}{2}(g(\nabla_X Y, Z) - g(\nabla_X JY, JZ)) = h_{XYZ}.$$

2.2°) $$h^2_{\overline{X}^h\overline{Y}^h\overline{Z}^v} = \frac{1}{4}(<[\overline{X}^h, \overline{Y}^h], \overline{Z}^v > + < [\overline{Z}^v, \overline{X}^h], \overline{Y}^h > +$$
$$+ < [\overline{Z}^v, \overline{Y}^h], \overline{X}^h > - < [\overline{X}^h, J_2\overline{Y}^h], J_2\overline{Z}^v > - < [J_2\overline{Z}^v, \overline{X}^h], J_2\overline{Y}^h > -$$
$$- < [J_2\overline{Z}^v, J_2\overline{Y}^h], \overline{X}^h >) = \frac{1}{4}(g(R(X,Y)U, Z) + g(R(X,JY)U, JZ)) =$$
$$= -\frac{1}{4}(g(R(X,Y)Z, U) + g(R(X,JY)JZ, U)).$$

By similar arguments we obtain

3.2°) $h^2_{\overline{X}^h\overline{Y}^v\overline{Z}^h} = -\frac{1}{4}(g(R(X,Z)Y, U) + g(R(X,JZ)JY, U)).$

4.2°) $h^2_{\overline{X}^v\overline{Y}^h\overline{Z}^h} = -\frac{1}{4}(g(R(Z,Y)X, U) - g(R(JZ,JY)X, U)).$

5.2°) $h^2_{\overline{X}^v\overline{Y}^v\overline{Z}^v} = 0$.

6.2°) $h^2_{\overline{X}^v\overline{Y}^v\overline{Z}^h} = 0.$

7.2°) $h^2_{\overline{X}^v\overline{Y}^h\overline{Z}^v} = 0.$

8.2°) $h^2_{\overline{X}^h\overline{Y}^v\overline{Z}^v} = \frac{1}{2}(g(\nabla_X Y, Z) - g(\nabla_X JY, JZ)) = h_{XYZ}.$

Here $h$ is the second fundamental tensor field of the pair $(J, g)$ on $M$.

## 3. Embeddings of almost Hermitian manifolds in almost hyperHermitian those

For an almost Hermitian manifold (*M*, *J*, *g*) we have constructed in **2** ahHs $(J_1, J_2, J_3, \hat{g})$ on *TM*. The manifold *M* can be considered as the null section $O_M$ in *TM* $(p \leftrightarrow o_p \in O_M \subset TM)$ and it is clear from (2.8) that $\hat{g}_{|M} = g$. All the results of **1** can be applied to a submanifold *M* in $(TM, \hat{g})$, see **[7]**. So, we can consider the normal tubular neighborhoods $Tb\left(M, \frac{\varepsilon(p)}{2}\right) \subset Tb(M, \varepsilon(p)) \subset TM$ and the deformations $\overline{J}_1, \overline{J}_2, \overline{J}_3, \overline{g}$ of the tensor fields $J_1, J_2, J_3, \hat{g}$ respectively.

**Theorem 2**. *Let (M, J, g) be an almost Hermitian manifold and* $Tb(M, \varepsilon(p))$ *be the corresponding normal tubular neighborhood with respect to* $\hat{g} = <\ ,\ >$ *on TM. Then M(*$O_M$*) is a totally geodesic submanifold of the almost hyperHermitian manifold* $\left(Tb\left(M, \frac{\varepsilon(p)}{2}\right), \overline{J}_1, \overline{J}_2, \overline{J}_3, \overline{g}\right)$, *where the ahHs* $(\overline{J}_1, \overline{J}_2, \overline{J}_3, \overline{g})$ *is the deformation of the structure* $(\overline{J}_1, \overline{J}_2, \overline{J}_3, \hat{g})$ *obtained in* **2°, 1.** *The structure* $(\overline{J}_1, \overline{g})$ *is Kaehlerian one.*

**Proof.** It follows from *theorem 1* that *M* is a totally geodesic submanifold of the Riemannian manifold $\left(Tb\left(M, \frac{\varepsilon(p)}{2}\right), \overline{g}\right)$.

Let $\widetilde{W}_0$ be a coordinate neighborhood in *TM* considered in **1°, 1.** A point $x \in \widetilde{W}_0$ has the coordinates $x_1, \ldots, x_n, x_{n+1}, \ldots, x_{2n}$ where $x_1, \ldots, x_n$ are coordinates of the point *p* in $\widetilde{V}_0 \subset M$ and $x_{n+1}, \ldots, x_{2n}$ are normal coordinates of *x* in $D\left(p, \frac{\varepsilon(p)}{2}\right)$.

We denote $X_i = \frac{\partial}{\partial x_i}$, $i = \overline{1,2n}$, $\hat{\nabla}_{X_i} X_j = \sum_k \hat{\Gamma}^k_{ij} X_k$, $\overline{\nabla}_{X_i} X_j = \sum_k \overline{\Gamma}^k_{ij} X_k$, $J X_j = \sum_k J^k_j X_k$, $\overline{J} X_j = \sum_k \overline{J}^k_j X_k$, $\hat{g}_{ij} = \hat{g}(X_i, X_j)$, $\overline{g}_{ij} = \overline{g}(X_i, X_j)$ where $\hat{\nabla}$ and $\overline{\nabla}$ are Riemannian connections of metrics $\hat{g}$ and $\overline{g}$, *J* is any tensor field from $J_1, J_2, J_3$.

Using the construction in **2°, 1** we have $\overline{g}_{ij}(x) = \hat{g}_{ij}(p)$, $\overline{J}^i_j(x) = J^i_j(p)$ on $Tb\left(M, \frac{\varepsilon(p)}{2}\right) \cap \widetilde{W}_0$. According to **[8]** we can write

$$(3.1) \qquad \sum_l \overline{g}_{lk} \overline{\Gamma}^l_{ij} = \frac{1}{2}\left(\frac{\partial \overline{g}_{kj}}{\partial x_i} + \frac{\partial \overline{g}_{ik}}{\partial x_j} - \frac{\partial \overline{g}_{ij}}{\partial x_k}\right)$$

It follows from (3.1) that $\overline{\Gamma}_{ij}^{l}(x)=\overline{\Gamma}_{ij}^{l}(p)$ and $\overline{\Gamma}_{ij}^{l}(x)=0$ i.e. $\overline{\nabla}_{X_i}X_j=0$ for $i=\overline{n+1,2n}$. Further, we get

$$\left(\overline{\nabla}_{X_i}\overline{J}\right)X_j=\overline{\nabla}_{X_i}\overline{J}X_j-\overline{J}\overline{\nabla}_{X_i}X_j=\sum_k\overline{\nabla}_{X_i}\overline{J}_j^kX_k-$$
$$-\overline{J}\left(\sum_k\overline{\Gamma}_{ij}^kX_k\right)=\sum_k\left(\overline{J}_j^k\overline{\nabla}_{X_i}X_k+\left(X_i\overline{J}_j^k\right)X_k\right)-$$
$$-\sum_{k,l}\overline{\Gamma}_{ij}^l\overline{J}_l^kX_k=\sum_{k,l}\left(\overline{J}_j^l\overline{\Gamma}_{il}^k-\overline{\Gamma}_{ij}^l\overline{J}_l^k+X_i\overline{J}_j^k\right)X_k,$$
$$\left(\left(\overline{\nabla}_{X_i}\overline{J}\right)X_j\right)(x)=\sum_{k,l}\left(\overline{J}_j^l\overline{\Gamma}_{il}^k-\overline{\Gamma}_{ij}^l\overline{J}_l^k+X_i\overline{J}_j^k\right)(x)X_{k|x}=$$
$$=\sum_{k,l}\left(\left(\overline{J}_j^l\overline{\Gamma}_{il}^k-\overline{\Gamma}_{ij}^l\overline{J}_l^k\right)(p)+\left(X_i\overline{J}_j^k\right)(x)\right)X_{k|x}.$$

It follows that $\overline{\nabla}_{X_i}\overline{J}=0$ for $i=\overline{n+1,2n}$.

For $i=\overline{1,n}$ $\left(X_i\overline{J}_j^k\right)(x)=\left(X_iJ_j^k\right)(p)$ and we obtain

$$\left(\left(\overline{\nabla}_{X_i}\overline{J}\right)X_j\right)(x)=\sum_{k,l}\left(J_j^l\hat{\Gamma}_{il}^k-\hat{\Gamma}_{ij}^lJ_l^k+X_iJ_j^k\right)(p)X_{k|x}.$$

From the other side we can write

$$\left(\left(\hat{\nabla}_{X_i}\overline{J}\right)X_j\right)(p)=\sum_{k,l}\left(J_j^l\hat{\Gamma}_{il}^k-\hat{\Gamma}_{ij}^lJ_l^k+X_iJ_j^k\right)(p)X_{k|p}.$$

According to **[3]** we have $\left(\overline{\nabla}_{X_i}J\right)X_j=\left(2h_{X_i}JX_j\right)(p)$ where the second fundamental tensor field *h* is defined by (2.11). From 1.1°) – 8.1°) it follows that $h_p^1=0$ for any $p\in M\left(U=o_p\in O_M\right)$. Thus, we have obtained $\overline{\nabla}J_1=0$ and the structure $\left(\overline{J}_1,\overline{g}\right)$ is Kaehlerian one on $Tb\left(M,\frac{\varepsilon(p)}{2}\right)$.

**QED.**

As a corollary we have got the following

**Theorem 3 [4]**. *Let (M, g) be a smooth Riemannian manifold and* $Tb(M,\varepsilon(p))$ *be the corresponding normal tubular neighborhood with respect to g = < , > on TM. Then M($O_M$) is a totally geodesic submanifold of the Kaehlerian manifold* $\left(Tb\left(M,\frac{\varepsilon(p)}{2}\right),\overline{J}_1,\overline{g}\right)$.

The classification given in **[5]** can be rewritten in terms of the second fundamental tensor field *h*, **[3].** Let *dimM* ≥ 6 and $2\beta(X)=\delta\Phi(JX)$, where $\Phi(X,Y)=g(JX,Y)$, then we have

| Class | Defining condition |
| --- | --- |
| **K** | $h = 0$ |
| **$U_1$ = NK** | $h_X X = 0$ |
| **$U_2$ = AK** | $\sigma h_{XYZ} = 0$ |
| **$U_3 = SK \cap H$** | $h_{XYZ} - h_{JXJYJZ} = \beta(Z) = 0$ |
| **$U_4$** | $h_{XYZ} = \frac{1}{2(n-1)}[<X,Y> \beta(Z) - <X,Z> \beta(Y) - <X,JY> \beta(JZ) + <X,JZ> \beta(JY)]$ |
| **$U_1 \oplus U_2$ = QK** | $h_{XYJZ} = h_{JXYZ}$ |
| **$U_3 \oplus U_4$ = H** | $N(J) = 0$ or $h_{XYJZ} = -h_{JXYZ}$ |
| **$U_1 \oplus U_3$** | $h_{XXY} - h_{JXJXY} = \beta(Z) = 0$ |
| **$U_2 \oplus U_4$** | $\sigma[h_{XYJZ} - \frac{1}{(n-1)} <JX,Y> \beta(Z)] = 0$ |
| **$U_1 \oplus U_4$** | $h_{XXY} = -\frac{1}{2(n-1)}[<X,Y> \beta(X) - \|X\|^2 \beta(Y) - <X,JY> \beta(JX)]$ |
| **$U_2 \oplus U_3$** | $\sigma[h_{XYJZ} + h_{JXYZ}] = \beta(Z) = 0$ |
| **$U_1 \oplus U_2 \oplus U_3$= = SK** | $\beta = 0$ |
| **$U_1 \oplus U_2 \oplus U_4$** | $h_{XYJZ} - h_{JXYZ} = \frac{1}{(n-1)}[<X,Y> \beta(JZ) - <X,Z> \beta(JY) + <X,JY> \beta(Z) - <X,JZ> \beta(Y)]$ |
| **$U_1 \oplus U_3 \oplus U_4$** | $h_{XJXY} + h_{JXXY} = 0$ |
| **$U_2 \oplus U_3 \oplus U_4$** | $\sigma[h_{XYJZ} + h_{JXYZ}] = 0$ |
| **U** | No condition |

**Proposition 4.** *Let (J, g) be from some class from the table above. Then the structure* $(\overline{J}_2, \overline{g})$ *has the analogous class on* $Tb\left(M, \frac{\varepsilon(p)}{2}\right)$.

**Proof.** From 1.2°) – 8.2°) it follows that $h^2_{\overline{XYZ}} = 2h_{XYZ}$. The rest is obvious from the table.

**4. Complex and hypercomplex numbers in differential geometry**

For the manifold $M$ we consider the products $M^2 = M \times M = \{(x; y) \mid x; y \in M\}$, $M^4 = M^2 \times M^2 = \{(x; y; u; v) \mid x; y, u; v \in M\}$ and the diagonals $\Delta(M^2) = \{(x; x) \in M^2\}$, $\Delta(M^4) = \{(x; x; x; x) \in M^4\}$. It is obvious that the manifold $\Delta(M^2)$ and $\Delta(M^4)$ are diffeomorphic to $M$ ($\Delta(M^2) \cong \Delta(M^4) \cong M$).

**Theorem 5 [6]**. *Let (M, ∇) be a manifold with a connection ∇ and $\pi : TM \to M$ be the canonical projection. Then there exists such a neighborhood $N_0$ of the null section $O_M$ in TM that the mapping*

$$\varphi : \pi \times \exp : X \to \left(\pi(X), \exp_{\pi(X)} X\right)$$

*is the diffeomorphic of $N_0$ on a neighborhood $N_\Delta$ of the diagonal $\Delta(M^2)$.*

Further, ∇ is a Riemannian connection of the Riemannian metric g. Combining the theorems 3, 5 we have obtained the following.

**Theorem 6**. *The diffeomorphism φ induces the Kaehlerian structure $(\overline{J}_1, \overline{g})$ on the neighborhood $N_\Delta$ of the diagonal $\Delta(M^2)$ and $\Delta(M^2) \cong M$ is a totally geodesic submanifold of the Kaehlerian manifold $(N_\Delta, \overline{J}_1, \overline{g})$.*

**Remark.** *Generally speaking, the complex structure of the Kaehlerian manifold $(N_\Delta, \overline{J}_1, \overline{g})$ is not compatible with the product structure of $M^2$. It means that if $z_l, l = \overline{1,n}$ are the complex coordinates of a point (x; y) $\in N_\Delta$, then, generally speaking, we can not find such real coordinates $x_l, y_l, l = \overline{1,n}$ of the points x, y ∈ M respectively that $z_l = x_l + iy_l$ where $i^2 = -1$.*

Combining the theorems 2, 3, 4, 5, 6 we have obtained the following.

**Theorem 7**. *There exists the hyperKaehlerian structure $(\overline{J}_1, \overline{J}_2, \overline{J}_3, \overline{g})$ on a neighborhood $\overline{N}_\Delta$ of the diagonal $\Delta(M^4)$* and *$\Delta(M^4) \cong M$ is a totally geodesic submanifold of the hyperKaehlerian manifold $(N_\Delta, \overline{J}_1, \overline{J}_2, \overline{J}_3, \overline{g})$.*

**Remark.** *Generally speaking, the hypercomplex structure of the hyperKaehlerian manifold $(\overline{N}_\Delta, \overline{J}_1, \overline{J}_2, \overline{J}_3, \overline{g})$ is not compatible with the product structure of $M^4$. It means that if $q_l, l = \overline{1,n}$ are the hypercomplex coordinates of a point (x; y; u; v) $\in \overline{N}_\Delta$, then, generally speaking we can not find such real coordinates $x_l, y_l, u_l, v_l$, $l = \overline{1,n}$ of the points x; y; u; v ∈ M respectively that $q_l = x_l + iy_l + ju_l + kv_l$ where $i^2 = j^2 = k^2 = -1$, $ij = -ji = k$.*

## References

1. S.A. Bogdanovich, A.A. Ermolitski, *On almost hyperHermitian structures on Riemannian manifolds and tangent bundles,* CEJM, 2(5) (2004) 615–623.
2. P. Dombrowski, *On the geometry of the tangent bundle,* J. Reine und Angew. Math., 210 (1962) 73–78.
3. A.A. Ermolitski, *Riemannian manifolds with geometric structures*, BSPU, Minsk, 1998 (in Russian), (English version: Internet, Google, Ermolitski )
4. A.A. Ermolitski, *Deformations of structures, embedding of a Riemannian manifold in a Kaehlerian one and geometric antigravitation*, Banach Center Publicantions, V. 76, Warszawa 2007, 505–514.
5. A. Gray, L. M. Herwella, *The sixteen classes of almost Hermitian manifolds and their linear invariants,* Ann. Mat. pura appl., 123 (1980) 35–58.
6. D. Gromoll, W. Klingenberg, W. Meyer, *Riemannsche geometrie im grossen*, Springer, Berlin, 1968 (in German).

**7.** M.W. Hirsch, *Differential topology. Graduate texts in mathematics*, 33, Springer, N.Y., 1976.
**8.** S. Kobayashi, K. Nomizu, *Foundations of differential geometry*, V. 1, Wiley, N.Y., 1963.
**9.** S. Kobayashi, K. Nomizu, *Foundations of differential geometry*, V. 2, Wiley, N.Y., 1969.